\documentclass[11pt]{article}
\title{On the moments of a polynomial in one variable}
\author{Michael M\"uger\footnote{Radboud University, Nijmegen, The Netherlands, {\tt mueger@math.ru.nl}}, Lars
Tuset\footnote{Oslo Metropolitan University, Oslo, Norway, {\tt larst@oslomet.no}}}

\usepackage{amssymb,amsmath, theorem,latexsym,ifthen,comment}
\usepackage{stmaryrd}
\usepackage{mathrsfs}
\usepackage{epic}
\usepackage{eepicemu}
\usepackage{epsfig}

\newlength{\dinwidth}
\newlength{\dinmargin}
\def\1#1{{\bf #1}}
\def\2#1{{\cal #1}}
\def\3#1{{\sl #1}}
\def\4#1{{\tt #1}}
\def\5#1{{\sf #1}}
\def\6#1{{\mathfrak #1}}
\def\7#1{{\mathbb #1}}
\def\8#1{{\mathscr #1}}

\newcommand{\be}{\begin{equation}}
\newcommand{\ee}{\end{equation}}
\newcommand{\ba}{\begin{array}}
\newcommand{\ea}{\end{array}}
\newcommand{\bea}{\begin{eqnarray}}
\newcommand{\eea}{\end{eqnarray}}
 \newcommand{\bean}{\begin{eqnarray*}}
\newcommand{\eean}{\end{eqnarray*}}

\newcommand{\ve}{\varepsilon}

\newcommand{\rarr}{\rightarrow}

\newcommand{\qed}{\hfill$\blacksquare$\\}

\def\endexem{\hfill{$\Box$}\medskip}

\theoremstyle{change}
\theoremheaderfont{\scshape}
\newtheorem{defin}{Definition}[section]
\newtheorem{defprop}{Definition/Proposition}[section]
\newtheorem{lemma}[defin]{Lemma}
\newtheorem{prop}[defin]{Proposition}
\newtheorem{theorem}[defin]{Theorem}
\newtheorem{coro}[defin]{Corollary}
\newtheorem{conj}[defin]{Conjecture}
\newtheorem{wish}[defin]{Desideratum}
\newtheorem{question}[defin]{Question}
\theorembodyfont{\rmfamily}
\newtheorem{remark}[defin]{Remark}
\newtheorem{exercise}[defin]{Exercise}

\newcommand{\bdefin}{\begin{defin}}
\newcommand{\blemma}{\begin{lemma}}
\newcommand{\bprop}{\begin{prop}}
\newcommand{\btheor}{\begin{theorem}}
\newcommand{\bcoro}{\begin{coro}}
\newcommand{\bdefprop}{\begin{defprop}}
\newcommand{\bexer}{\begin{exercise}}
\newcommand{\edefprop}{\end{defprop}}
\newcommand{\edefin}{\end{defin}}
\newcommand{\elemma}{\end{lemma}}
\newcommand{\eprop}{\end{prop}}
\newcommand{\etheor}{\end{theorem}}
\newcommand{\ecoro}{\end{coro}}
\newcommand{\eexer}{\end{exercise}}
\newcommand{\bconj}{\begin{conj}}
\newcommand{\econj}{\end{conj}}
\newcommand{\bwish}{\begin{wish}}
\newcommand{\ewish}{\end{wish}}
\newcommand{\bquestion}{\begin{question}}
\newcommand{\equestion}{\end{question}}
\newcommand{\brem}{\begin{remark}}
\newcommand{\erem}{\endexem\end{remark}}

\newcommand{\prf}{{\noindent\it Proof. }}
\def\mobj#1{\raise .4\unitlens\hbox{\put(0,0){$#1$}}}
\def\mychi{\raise 2pt\hbox{$\chi$}}

\begin{document}

\maketitle

\abstract{Let $f$ be a non-zero polynomial with complex coefficients and define
$M_n(f)=\int_0^1f(x)^n\,dx$. We use ideas of Duistermaat and van der Kallen to prove
$\limsup_{n\rarr\infty}|M_n(f)|^{1/n}>0$. In particular, $M_n(f)\ne 0$ for infinitely many
$n\in\7N$.}

%%%%%%%%%%%%%%%%%%%%%%%%%%%%

\section{Continuous real valued functions}
While our main interest is in complex polynomials, cf.\ the Abstract, it is instructive to first
have a quick look at the much simpler case of real valued functions.

\bprop Let $I=[0,1]$ and $f\in C(I,\7R)$. For $n\in\7N$ define
\be M_n(f) = \int_0^1 f(x)^n dx, \quad\quad\quad M(f) = \sup_I|f|. \label{eq-M}\ee
(We will often just write $M_n, M$.) Then
\[ \limsup_{n\rarr\infty} |M_n(f)|^{1/n}=M(f).\]
\eprop

\prf We first assume $f(I)\subseteq[0,\infty)$. By compactness there is an $x_0\in I$ such that 
$f(x_0)=M$. Let $\ve>0$. By continuity there is $\delta>0$ such that  
$f(x)\ge M-\ve$ on $(x_0-\delta,x_0+\delta)$. We may assume $\delta$ small enough so that at least
one of the intervals $(x_0-\delta,x_0],\,[x_0,x_0+\delta)$ is contained in $I$. In view of $f\ge 0$,
we have $\int_I f^n\ge(M-\ve)^n\delta\ \forall n$. This implies $\lim(\int f^n)^{1/n}\ge M-\ve$.
Since $\ve>0$ was arbitrary, the limit is $\ge M$. The converse inequality being trivial, we have 
$(\int f^n)^{1/n}\rarr M$.

For general $\7R$-valued $f$, applying the above reasoning to $g=f^2$ gives
\[ \lim_{n\rarr\infty} M_{2n}(f)^{1/n}=\lim_{n\rarr\infty} \left(\int_I g^n \right)^{1/n}=M^2,
\]
which together with $|M_n|\le M^n$ proves the claim.
\qed

No result of comparable generality seems to be known for continuous complex valued functions. We
therefore now turn to polynomials.

\begin{comment}
\brem The question whether $\lim_{n\rarr\infty}|M_n(f)|^{1/n}$ exists (then clearly equaling
$M(f)$) is harder. Sufficient conditions are $M_+\ne|M_-|$ and ($M_+=|M_-|$, but
$\lambda(f^{-1}(M))\ne\lambda(f^{-1}(-M))$). The function $f(x)=x-1/2$ satisfies neither of these 
conditions, and in fact $\lim_{n\rarr\infty}|M_n(f)|^{1/n}$ does not exist since $M_n(f)=0$ for all
odd $n$.
\erem
\end{comment}

\begin{comment}
\prf Assume $M_+\ne |M_-|$. In view of $f_+(x)f_-(x)=0$ for all $x$, it is clear that
$M_n(f)=M_n(f_+)+M_n(f_-)$. Assume $M_+>|M_-|$, and let $\ve=(M_+-|M_-|)/2$. Then as in the proof of
(i) there is a $\delta>0$ such that $M_n(f_+)\ge (M_+-\ve)^n \delta$, while
$|M_n(f_-)|\le|M_-|^n$. In view of our choice of $\ve$, we have
$\frac{M_n(f_-)}{M_n(f_+)}\rarr 0$. Thus the contribution of $f_+$ to $M_n(f)$ drowns out that  
of $f_-$, and taking $n$-th roots the result follows. Analogously if $|M_-|>M_+$.

Now assume $M_+=|M_-|$, but $\lambda(f^{-1}(M))\ne\lambda(f^{-1}(-M))$. Then
\bean M_n(f) &=& M^n \mu_++(-M)^n\mu_-+\int_{|f|^{-1}([0,M))} f^n(x)dx \\
  &=& M^n\left( \mu_+ +(-1)^n\mu_- +\int_{|f|^{-1}([0,M))} \left(\frac{f(x)}{M}\right)^n dx\right)
\eean
The integral tends to zero by the dominated convergence theorem. Thus
\[ |M_n(f)|^{1/n}=M \left| \mu_+ +(-1)^n\mu_- + o(1)\right|^{1/n}, \]
which converges to $M$ since $|\mu_+ +(-1)^n\mu_-|\ge |\mu_+-\mu_-|>0$.
\qed

\end{comment}

%%%%%%%%%%%%%%%%%%%%%%%%%%%%%%%%%%%%%%%%%%%%%%%%%%%%%%%%%%%%%%%%%%%%%%%%%%%%%%%%%%%%%%%%%%%%%%%%%%%%%%%%%%%%%%%%%%

\section{Complex polynomials}
In this section, we consider polynomials with complex coefficients on $I=[0,1]$. We define $M_n(f)$
and $M(f)$ as in (\ref{eq-M}). Our aim is to prove:

\btheor \label{theor-X} If $f$ is a non-zero polynomial then
$\limsup_{n\rarr\infty}|M_n(f)|^{1/n}>0$. In particular, $M_n(f)\ne 0$ for infinitely many $n\in\7N$.
\etheor

\brem The second statement was already proven in \cite[Corollary 4.1]{FPYZ} by mostly algebraic
methods, involving some amount of number theory (including Dirichlet's theorem on primes in
arithmetic progressions!).
That result recently found an application \cite{dings} in a tentative approach to the Mathieu
conjecture \cite{mathieu} for $SU(2)$. 
Our proof of Theorem \ref{theor-X} will be purely analytic, inspired by
the proof of the Mathieu conjecture for $S^1$ given in \cite[Theorem 2]{dvdk}.
\erem

\prf We first observe that if $f=C=$\,const then $M_n=C^n\ \forall n$, which clearly implies the
theorem in this case. Thus from now on we may and will assume $d=\deg(f)\ge 1$, thus in particular
$f\ne 0$. It will suffice to prove the theorem for monic polynomials since $M_n(cf)=c^nM_n(f)$.

We will first prove that at least one moment is non-zero.
In view of $|M_n|\le M^n$, the generating function $F(t)=\sum_{n=0}^\infty t^n M_n$ clearly converges
to a holomorphic function on the open disc $B_{1/M}(0)$, and on this domain
\be F(t)=\int_0^1 \frac{dx}{1-tf(x)}. \label{eq-F1}\ee
We define a finite subset of $\7C$ by
\be S=\{ f(z)\ | \ z\in\7C, f'(z)=0\}\cup\{ f(0), f(1)\}.\label{eq-S}\ee
We will prove that $F(t)$ can be analytically continued along any path from small $t$ (where $F(t)$
is known to be analytic) to $\infty$ such that $\tau=1/t$ avoids the set $S$. We will then show that
$F(t)\rarr 0$ along any such path. This will show that $F\ne 1$, so that at least one moment
$M_n(f)$ is non-zero. In view of
\be F(t)=\frac{-1}{t}\int_0^1 \frac{dx}{f(x)-\tau}, \label{eq-F2}\ee
we need some information about the zeros of the polynomial $f-\tau$.

Since we assume $d=\deg(f)\ge 1$, for each $\tau\in\7C$ the equation $f(z)=\tau$ has $d$ solutions
$z_{1,t},\ldots, z_{d,t}$. As long as these solutions are pairwise distinct, they depend
holomorphically on $\tau$. (This is 
usually proven using Rouch\'e's theorem, but it can also been construed as an application of the
holomorphic implicit function theorem, see e.g.\ \cite[Section 6.1]{krantz}.) The condition that all
zeros of $f-\tau$ be distinct is equivalent to none of them being multiple, which is to say none of
them is a critical point of $f-\tau$ or, equivalently, of $f$. This in turn is equivalent to $\tau$
not being a critical value of $f$. Thus as $\tau$ traces a path in $\7C\backslash S$, the functions
$z_{1,t},\ldots, z_{d,t}$ are analytic in a neighborhood of the path. (These functions are
multivalued in the sense that they depend on the homotopy class of the chosen path, but we choose
one path once and for all.)

Since we insist on $\tau\not\in S$, the zeros of the monic polynomial $f-\tau$ are pairwise
distinct, so that we can apply partial fraction expansion in its most basic form: 
\[ \frac{1}{f(x)-\tau}=\sum_{k=1}^d \left(\prod_{\ell\ne k}\frac{1}{z_{k,t}-z_{\ell,t}}\right)\frac{1}{x-z_{k,t}}. \] 
Formally integrating over $[0,1]$ gives 
\be F(t)=\frac{-1}{t}\sum_{k=1}^d \left(\prod_{\ell\ne k}\frac{1}{z_{k,t}-z_{\ell,t}}\right) \log\frac{1-z_{k,t}}{-z_{k,t}}.
\label{eq-F3}\ee
(We will return to the choice of branches for the logarithms shortly.) 
As a consequence of $\tau\not\in S$, the denominators $z_{k,t}-z_{\ell,t}$ never vanish (for finite
$t$). Since the forbidden set $S$ also contains the endpoint values $f(0)$ and $f(1)$, the
assumption $\tau\not\in S$ implies $f(0)\ne\tau\ne f(1)$. Since the $z_{i,t}$ satisfy
$f(z_{i,t})=\tau$, we conclude that as long as $\tau\not\in S$, the solutions
$z_{1,t},\ldots,z_{d,t}$ of $f(z)=\tau$ assume neither of the values $0$ or $1$. Thus the arguments
of the logarithms in (\ref{eq-F3}) are finite and non-zero. For small $t$, where $F$ is defined a
priori, $\tau=1/t$ is large, thus is not in $S$. Now the branches of the logarithms can be chosen such
that (\ref{eq-F3}) holds for $t\in B_{1/M}(0)$. As $t$ increases, the analytic continuation of the 
logarithms is done by lifting the path traced by $\frac{1-z_{k,t}}{-z_{k,t}}$ to the Riemann surface
of the logarithm.

We have now achieved our first goal of analytically extending $F$ to a (multi-valued) analytic
function on $\7C\backslash S$. It remains to study the behavior of $F(t)$ as $t\rarr\infty$.

We first assume $0\not\in S$. Thus $f(0)\ne 0\ne f(1)$, and $0$ is not a critical value. The latter
is equivalent to $f$ having no repeated zeros. Under this assumption, the solutions $z_{i,t}$ of
$f(z_{i,t})=\tau$ extend holomorphically to $\tau=0$, i.e.\ $t=\infty$. The limits
$z_{1,\infty},\ldots, z_{d,\infty}$ are the zeros of $f$, thus they are all distinct by our 
assumption. The consequence $f(0)\ne 0\ne f(1)$ of the assumption 
$0\not\in S$ implies $z_{i,\infty}\not\in\{0,1\}$ for all $i$. Thus all terms in (\ref{eq-F3}) after
the $-1/t$ have finite limits as $t\rarr\infty$, so that $F(t)=O(1/t)$. As explained earlier, this
implies $F\ne 1$, thus $M_n(f)\ne 0$ for at least one $n$.

Now assume $0\in S$. This can arise from the existence of a multiple zero of $f$ or from the
vanishing of $f$ at $0$ or at $1$, or any combination of these. Even if this is the case, the solutions
$z_{i,t}$ of $f(z_{i,t})=\tau$ considered above do converge to zeros of $f$ as $\tau\rarr 0$. If $f$
has zeros with multiplicities, we must be more careful in studying the $t\rarr\infty$ limit of
(\ref{eq-F3}) since some of the denominators $z_{k,t}-z_{\ell,t}$ will tend to zero. And if $f(0)=0$
then we will have $z_{i,t}\rarr 0$ for some $i$ as $t\rarr\infty$, so that the behavior of the
corresponding logarithmic factor must be reconsidered. Similar considerations apply when $f(1)=0$.

Assume $f$ has an $n$-fold zero at $z=0$. Thus $f(z)=z^ng(z)$, where $g(0)\ne 0$. For small $z$ we
have $f(z)=z^n g(0)+O(z^{n+1})$. Now, precisely $n$ of the solutions $z_{i,t}$ of $f(z)=\tau$ will
tend to zero as $\tau\rarr 0$. For these $i$ we have $1/t=\tau=f(z_{i,t})\sim z_{i,t}^n g(0)$ for
small $z_{i,t}$, equivalently large $t$. Thus the $z_{i,t}$ behave like the $n$ $n$-th roots of
$\frac{1}{tg(0)}$ as $t\rarr\infty$. The conclusion important for us is that each of the functions
$z_{i,t}$ that tend to zero as $t\rarr\infty$, do so like $t^{-1/n}$ times a non-zero
constant. Plugging this into $\log\frac{1-z_{i,t}}{-z_{i,t}}$ gives a proportional to $\log t$ as
$t\rarr\infty$. (An entirely similar analysis arrives at the same conclusion in the case $f(1)=0$, so
that at least one $z_{i,t}$ tends to $1$ as $t\rarr\infty$.) Since $\frac{\log t}{t}\rarr 0$ as
$t\rarr\infty$, we will still be able to conclude that $F(t)\rarr 0$ as long as the denominators do
not create problems. Since this clearly cannot happen if $f$ does not have multiple zeros, we have
generalized our proof of $F(t)\rarr 0$ to this case, whether or not $f$ vanishes at zero or one.

Turning to the multiple zero case, we first observe that the $k$-summand in (\ref{eq-F3}) will
behave nicely, i.e.\ like $\frac{1}{t}$ or $\frac{\log t}{t}$, as $t\rarr\infty$, provided
$z_{k,\infty}$ is a simple zero of $f$, since then
$\prod_{\ell\ne k}(z_{k,t}-z_{\ell,t})\rarr\prod_{\ell\ne k}(z_{k,\infty}-z_{\ell,\infty})\ne 0$.

It remains to study the $k$-terms in (\ref{eq-F3})
for which $z_{k,t}$ converges to a multiple zero of $f$. Let $n\ge 2$ be the multiplicity of the
zero $z_{k,\infty}$. We know that the factor with the logarithm
will either converge to a non-zero constant  or behave like $\log t$, depending on whether
$z_{k,\infty}\not\in\{0,1\}$ or not. The factor $z_{k,t}-z_{\ell,t}$ will converge to 
$z_{k,\infty}-z_{\ell,\infty}$, which is non-zero unless also $z_{\ell,t}\rarr z_{k,\infty}$, which
of course is possible since we now allow $f$ to have multiple zeros. By the argument given above for
a multiple zero of $f$ at $0$, applied to $g(w)=f(z_0+w)$, we find: If $f$ has an $n$-fold zero at
$z_0$, then for the $z_{k,t}$ that tend to $z_0$ as $t\rarr\infty$ 
we have that $z_{k,t}-z_0$ behaves like a non-zero constant (the same for all $k$ under
consideration) multiplied by $t^{-1/n}$ and an $n$-th root of unity. The distance between distinct
$n$-th roots of unity (for fixed $n$) obviously is bounded below by a positive constant. Thus for 
the $z_{\ell,t}$ that also tend to $z_{k,\infty}$, the difference $z_{k,t}-z_{\ell,t}$ is bounded
below by $t^{-1/n}$ times a positive constant as $t\rarr\infty$. We conclude that the $\ell$-factor
in the $k$-summand in (\ref{eq-F3}) is bounded by a constant times $t^{1/n}$ as $t\rarr\infty$ for
each $\ell$ with $z_{k,\infty}=z_{\ell,\infty}$. There are precisely $n-1$ of these (since $k$
itself does not contribute to the product). Thus the overall behavior of the $k$-summand is
$t^{-1+(n-1)/n}$, modified by $\log t$ if $z_{k,\infty}\in\{0,1\}$. Since $-1+(n-1)/n=-1/n<0$, we
see that all $k$-summands in (\ref{eq-F3}) tend to zero as $t\rarr\infty$. This proves $F(t)\rarr 0$
in full generality.

In the preceding argument, the multiplicity $n$ was of course bounded by the degree $d$ of $f$. Thus
$F(t)=O(t^{-1/d}\log t)$, which holds uniformly in the direction once
$|t|>\min(\{ |s|\ | \ s\in S\backslash\{0\}\})^{-1}$. Assume now 
$\limsup_{n\rarr\infty}|M_n(f)|^{1/n}=0$. This would imply that the power series defining $F$ has 
infinite convergence radius so that $F$ is an entire function. Together with the uniform
decay at infinity this clearly means that $F$ is a bounded entire function and therefore constant by
Liouville's theorem. But this is inconsistent with $F(0)=1$ and $F(t)\rarr 0$ at infinity. This
contradiction proves $\limsup_{n\rarr\infty}|M_n(f)|^{1/n}>0$.
\qed

The fact that the power series defining $F$ has finite radius of convergence $R$ implies that the
function $F$ must have a singularity at some $t\in\7C$ with $|t|=R$. In view of (\ref{eq-F3}), such
a singularity can occur only when $t=1/s$ for a non-zero $s\in S$. This implies
\[ \limsup_{n\rarr\infty}|M_n(f)|^{1/n}\le \max\{ |s| \ | \ s\in S\}. \]
It would clearly be desirable to have more precise results about
$\limsup_{n\rarr\infty}|M_n(f)|^{1/n}$. 
In some cases it it easy to show that equality occurs in the above inequality. This definitely hapens
when $|f(0)|$ is larger than the absolute values of the other elements of $S$, or when
$f(0)\ne f(1)$ and $|f(0)|=|f(1)|$ is larger than the absolute values of the critical values. In
these cases the singularity arising in (\ref{eq-F3}) from the
vanishing of the argument of the logarithm cannot be offset by
something else. One might conjecture the following, but without too
much confidence: 

\bconj Let $f$ be a polynomial and define $S$ as in (\ref{eq-S}). Then
\[ \limsup_{n\rarr\infty}|M_n(f)|^{1/n}=\max\{ |s| \ | \ s\in S\}. \]
\econj

\vspace{.5cm}

\noindent{\it Acknowledgment.} M.\ M.\ thanks L.\ T.\ and OsloMet for hospitality and financial support.


\begin{thebibliography}{99}
\bibitem{dings} T. Dings, E. Koelink: {\it On the Mathieu conjecture for $SU(2)$}.
Indag. Math. (N.S.) {\bf 26}, 219-224 (2015).
\bibitem{dvdk} J. J. Duistermaat, W. van der Kallen: Constant terms in powers of a Laurent
polynomial. Indag. Math. (N.S.) {\bf 9}, 221-231 (1998).
\bibitem{FPYZ} J.P. Francoise, F. Pakovich, Y. Yomdin, W. Zhao: Moment vanishing problem and positivity: Some examples.
Bull. Sci. math. {\bf 135}, 10-32 (2011).
\bibitem{krantz} S. G. Krantz, H. R. Parks: {\it The Implicit Function Theorem. History, theory, and
applications}. Birkh\"auser, 2003.
\bibitem{mathieu} O. Mathieu: Some conjectures about invariant theory and their
applications. pp. 263-279 {\it in:} J. Alev, G. Cauchon (eds.): {\it Alg\`ebre non commutative,
groupes quantiques et invariants.} (Proceedings of the 7th Franco-Belgian Conference, Reims, June
26-30, 1995.). Soc. Math. France, 1997. 
\end{thebibliography}
\end{document}